\documentclass[a4paper,11pt]{amsart}
\usepackage{hyperref,latexsym}
\usepackage{enumerate}

\theoremstyle{plain}
\newtheorem{theorem}{Theorem}[section]

\theoremstyle{definition}

\newtheorem{example}{Example}
\theoremstyle{remark}
\newtheorem{remark}{Remark}

\begin{document}

\title[Algebraic Birkhoff conjecture]
      { Algebraic Birkhoff conjecture for billiards on Sphere and Hyperbolic plane }

\date{1 March 2016}
\author{Misha Bialy and Andrey E. Mironov}
\address{M. Bialy, School of Mathematical Sciences, Raymond and Beverly Sackler Faculty of Exact Sciences, Tel Aviv University,
Israel} \email{bialy@post.tau.ac.il}
\address{A.E. Mironov, Sobolev Institute of Mathematics and
 } \email{mironov@math.nsc.ru}
\thanks{M.B. was supported in part by ISF grant 162/15 and A.E.M. was
supported by RSF (grant 14-11-00441). It is our pleasure to thank
these funds for the support}

\subjclass[2000]{ } 

\begin{abstract}We consider a convex curve
$\gamma$ lying on the Sphere or Hyperbolic plane. We study the
problem of existence of polynomial in velocities integrals for
Birkhoff billiard inside the domain bounded by $\gamma$. We extend
the result by S. Bolotin (1992) and get new obstructions on
polynomial integrability in terms of the dual curve $\Gamma$. We
follow a method which was introduced by S. Tabachnikov for Outer
billiards in the plane and was applied later on in our recent paper
to Birkhoff billiards with the help of a new the so called Angular
billiard.

\end{abstract}

\maketitle

\section{Motivation and the results}
It is an old and difficult problem in theory of billiards going back
to Birkhoff to describe those billiard tables $\Omega$ which admit a
nonconstant function on the unit tangent bundle $T_1\Omega$ which is
polynomial in velocities and keeps constant values along the orbits
of the billiard flow. Traditionally such a function is called
\textit{Polynomial Integral}. Under certain condition on the
boundary curve, a remarkable result on polynomially integrable
billiards in the Euclidean plane was obtained by Bolotin (1990) in
\cite{B1}.  Later on, it was extended \cite{B2} also for billiards
on constant curvature surfaces.

For Outer billiards the problem about polynomial integrals in the
plane was recently considered \cite{T}. In our recent paper
\cite{BM}, the new model of Angular billiard was introduced. Using
Angular billiard we extended the method by S. Tabachnikov for
Birkhoff billiards in the Euclidean plane. Obstructions to
integrability which we obtained by this method are in a sense dual
to those obtained by Bolotin.

In this, paper we consider Birkhoff billiard inside a convex domain
$\Omega$ lying on the surface $\Sigma$ of constant curvature $\pm
1,$ i.e on standard Sphere or Hyperbolic plane. This is a particular
case of Riemannian billiard. Point moves along geodesic inside
$\Omega$ and hitting the boundary reflects according to the low of
geometric optics. It turns out, that for the Sphere and Hyperbolic
plane it is possible to follow our approach of \cite{BM} in order to
find the obstructions to polynomial integrability in terms of the
dual curve of the boundary, extending the result of \cite{B2} for
the constant curvature case.

In what follows, for the case of $K=1$ we realize $\Sigma$ as the
standard unit sphere with the induced metric from Euclidean
$\mathbb{R}^3$, and as the upper sheet of the hyperboloid
$\{x_1^2+x_2^2-x_3^2=-1\}$ in $\mathbb{R}^3$ endowed with
$ds^2=dx_1^2+dx_2^2-dx_3^2$, for the case $K=-1.$ Let $\gamma$ be a
smooth regular arc of the boundary of $\Omega.$ We shall denote by
$\hat{\gamma}$ the image of the curve $\gamma$ under the projection
$\mathbb{R}^3\setminus\{0\}\rightarrow\mathbb{R}P^2$.
\begin{example}
\label{example} Consider as an example the domain $\Omega$ inside
the ellipse on $\Sigma$. In this case the boundary curve $\gamma$ is
the intersection of $\Sigma$ with a quadratic cone
$\{ax^2+by^2+cz^2=0\}$. This cone defines a curve $\hat{\gamma}$ in
the projective space. Birkhoff billiard inside $\Omega$ is known to
admit an additional  integral which is quadratic in velocities (see
\cite{B2}, \cite{V}).
\end{example}

Let us formulate first the key result of \cite{B2}, needed for the
classification of billiards on $\Sigma$ admitting polynomial
integrals. Everywhere in the paper we shall denote by
$\Lambda\subset\mathbb{C}P^2$ the absolute, defined by:
$$\Lambda=\{(x_1:x_2:x_3):x_1^2+x_2^2\pm x_3^2=0\},$$
where $+$ sign is taken for the Sphere and $-$ sign for the
Hyperboloid.

\begin{theorem} (Bolotin, 1992) {\it Let $\gamma$ be a smooth non-geodesic arc of the
boundary curve of the domain $\Omega\subset\Sigma$. Suppose that
Birkhoff billiard inside $\Omega\subset\Sigma$ admits a non-constant
polynomial integral $\Phi$ on the energy level $\{|v|=1\}.$ It then
follows that
 $\hat{\gamma}$ is necessarily algebraic curve.
Moreover, let $\tilde{\gamma}$ in $\mathbb{C}P^2$
 be the irreducible component of $\hat{\gamma}.$
 If  $\tilde{\gamma}$ is a smooth curve, such that at least one intersection point of $\tilde{\gamma}$ with the absolute $\Lambda$ is transversal,
 then $\tilde{\gamma}$ is of degree 2.
 }
\end{theorem}
It turns out to be a very difficult and open question how to remove
the severe assumption in Bolotin's theorem. For instance, one would
expect that the Example \ref{example}  exhausts all convex domains
with smooth boundaries having additional polynomial integral. We
refer to this problem as Algebraic Birkhoff conjecture.

Our approach gives the following result in this direction:
\begin{theorem}\label{main}Let $\tilde{\Gamma}$ in $\mathbb{C}P^2$ be the dual curve of
$\tilde{\gamma}$. Then the following alternative holds: either
$\tilde{\Gamma}$ is of degree 2, or $\tilde{\Gamma}$ necessarily
contains singular points, so that all the singular points and all
inflection points of $\tilde{\Gamma}$ belong to the absolute
$\Lambda$.
\end{theorem}

This result is motivated by our approach in \cite{BM} for the planar
case, where we introduced the so called Angular billiard---the
object dual to the Birkhoff billiard near the boundary curve. Also
in the present paper the proof of Theorem \ref{main} is based on an
identity (\ref{pme}) obtained in terms of the dual billiard. Then we
extract first nontrivial term of the power series expansion of the
identity in the parameter $\varepsilon$ and get a remarkable
equation (\ref{Hess}) on $\tilde{\Gamma}$. The latter is studied in
algebro-geometric way.

Let us denote by $\Gamma$ the dual curve to $\gamma$. Let us recall
that the duality in geometric terms can be defined as follows. For a
geodesic on $\Sigma$ with a velocity vector $v\in T_r\Sigma$
corresponds the point $M=r\wedge v$, i.e. the momentum of $v$ (here
and later $\wedge$ stands for the standard vector product in
$\mathbb{R}^3$). Since the momentum is preserved along geodesics,
then this correspondence is well defined.  Moreover, one can easily
see that it extends to the usual projective duality.

By the definition, in the case of the sphere the dual point also
lies on the sphere and for the Hyperboloid the dual point belongs to
de Sitter surface, i.e. to one sheeted Hyperboloid
$\{x_1^2+x_2^2-x_3^2=1\}$. Choosing positive orientation and the
parametrization by arc-length on $\gamma$ we get a parametrization
of the dual curve:
$$\Gamma(s)= \gamma(s)\wedge\dot{\gamma}(s).$$

In both cases our approach is based on the following
\begin{theorem}
\label{pm}Assume that Birkhoff billiard in $\Omega$ admits a
polynomial integral $\Phi$ of even degree $n$. Then there exists a
homogeneous polynomial $\Psi$ of degree $n$ vanishing on $\Gamma$
which satisfies the the following identity. For any point $M\in
\Gamma$ and a non zero tangent vector $w\in T_{M}\Gamma$ one has:
\begin{equation}\label{pme}
\Psi(M-\varepsilon w)=\Psi(M+\varepsilon w),
\end{equation}
for all sufficiently small $\varepsilon.$
\end{theorem}
The geometric sense of this result can be expressed by saying that
the dual object to Birkhoff billiard on constant curvature surface
is the Outer billiard \cite{Tbook}. In the case of the sphere this
is well known to billiard players community but slightly less
obvious for the case of Hyperboloid. We shall make this precise in
the remark after the proof of Theorem \ref{pm}

The algebraic part of our approach is based on the following
theorem.

Let us denote by $F\in \mathbb{C}[x_1,x_2,x_3]$ the irreducible
homogenous defining polynomial of the component
$\tilde{\Gamma}\subset\mathbb{C}P^2$, and let $d=\textrm{deg} F,\
d\geq 2$. As usual we introduce
$$
 {\rm Hess}(F)(x_1,x_2,x_3):=\det\left(
  \begin{array}{ccc}
    F_{x_1x_1} & F_{x_1x_2} & F_{x_1x_3} \\
    F_{x_2x_1} & F_{x_2x_2} & F_{x_2x_3} \\
    F_{x_3x_1} & F_{x_3x_2} & F_{x_3x_3}\\
  \end{array}
\right).
$$
\begin{theorem}\label{algebra}Suppose that Birkhoff billiard in $\Omega$ admits a
polynomial integral $\Phi$ of even degree $n$. Let $F$ be the
irreducible homogeneous defining polynomial of $\tilde{\Gamma},
d=\rm {deg} F$. Then the following identity holds:
\begin{equation}\label{Hess}
Q^3({(\rm Hess(F))}^k-c(x_1^2+x_2^2+Kx_3^2)^{\alpha}=F\cdot R,\quad
\forall(x_1,x_2,x_3)\in\mathbb{C}^3.
\end{equation}
Here $\alpha$ is a positive integer, when $d>2$; $K$ is the
curvature $\pm1$ and $Q,R$ are homogeneous polynomials, so that $Q$
does not vanish on $\tilde{\Gamma}$ identically.
\end{theorem}
Postponing the proof of Theorem \ref{pm} and \ref{algebra} let us
prove Theorem \ref{main}.
\begin{proof}
We follow the idea of Lemma 3 of \cite{T}. Consider the situation in
$\mathbb{C}P^2.$ Any intersection point in $\mathbb{C}P^2$ between
Hessian curve of ${\rm Hess}(\tilde{\Gamma})$ with $\tilde{\Gamma}$
is either singular or inflection point of $\tilde{\Gamma}$. So, if
there is a singular or inflection point $(a:b:c)\in\tilde{\Gamma}$
such that $a^2+b^2+Kc^2\ne 0,$ it then follows from (\ref{Hess})
that $c=0.$ Therefore, ${\rm Hess}(F)$  must vanish on
$\tilde{\Gamma}$ identically,  since $Q$ does not vanish identically
on $\tilde{\Gamma}.$ This implies that $\tilde{\Gamma}$ is a line
(see \cite{F}), but this is impossible.

Let us prove now that $\tilde{\Gamma}$ must have singular points. If
on the contrary $\tilde{\Gamma}$ is a smooth curve, then it follows
from (\ref{Hess}) that all inflection points must belong to the
absolute $\Lambda$ defined by the equations
$$
 \Lambda=\{ (x_1:x_2:x_3)\in\mathbb{C}P^2:{x_1^2+x_2^2+Kx_3^2=0} \}.
$$
Recall, $d$ is the degree of $\tilde{\Gamma}.$ Then the Hessian
curve intersects $\tilde{\Gamma}$ exactly in inflection points, and
moreover, it is remarkable fact that the intersection multiplicity
of such a point of intersection equals exactly the order of
inflection point (see \cite{W}), and hence does not exceed  $(d-2).$
Furthermore, the absolute $\Lambda$ intersects $\tilde{\Gamma}$
maximum in $2d$ points together. Hence,  we have altogether counted
with multiplicities not more than $ 2d(d-2),$ but on the other hand
the Hessian curve has degree $3(d-2)$ and thus by Bezout theorem the
number of intersection points with multiplicities is $3d(d-2).$ This
contradiction shows that $\tilde{\Gamma}$ can not be a smooth curve
unless $d=2.$ Theorem \ref{main} is proved.
\end{proof}
The plan of the rest of the paper is as follows: in Section
\ref{geometric} we give the geometric proof of Theorem \ref{pm} and
then in Section \ref{sphere} we treat both cases, of Sphere and
Hyperboloid in separate Subsections.
\section{Proof of Theorem\ref{pm}}
\label{geometric} In this section we prove Theorem \ref{pm}.

The first observation is the following \cite{B1}, \cite{B2},
\cite{KT}. Let $\Phi$ is a polynomial integral of Birkhoff billiard
inside $\Omega$. Then, one can assume, with no loss of generality,
that $\Phi(r,v)$ is homogeneous in velocities of certain even degree
$n$. Moreover, since components of the Momentum $M=r\wedge v$ are
preserved under the geodesic flow on $\Sigma$ one can show that
$\Phi(r,v)\equiv\Psi(M),$ where $\Psi$ is a homogeneous polynomial
in the components of $M$ of degree $n.$ Furthermore, by a simple
modification of $\Phi$ one can achieve that $\Phi$ vanishes on
tangent vectors to the boundary curve $\gamma$. This implies that
the polynomial $\Psi$ vanishes on the dual curve $\Gamma$, because
the latter was constructed by:
$$\Gamma(s):=M(s)=\gamma(s)\wedge\dot{\gamma}(s),$$
with the help of arc-length parameter $s$ on $\gamma$. We remind
that Riemannian metric on $\Sigma$ is induced from the Euclidean
Space in the case of Sphere, and from Minkowski metric in the case
of Hyperboloid.

Fix a point $r=\gamma(s)$ and fix the orthonormal frame $\{v,n \}$
in $T_{r}\Sigma$,  $v=\dot{\gamma}(s)$ and $n$ is the positive unite
normal. Denote the dual point by $$M=\Gamma(s)=\gamma(s)\wedge
v=r\wedge v.$$

Let us compute the tangent vector to the dual curve at $M$:
\begin{equation}\label{N}
\dot{\Gamma}(s)=\dot{M}(s)=\gamma(s)\wedge\ddot{\gamma}(s)=
\gamma\wedge(\nabla_v\dot{\gamma}+N),
\end{equation}
where $\nabla$ is the Levi-$\breve{\rm C}$ivita connection of the
induced metric on the surface $\Sigma$, $N$ is a normal vector to
$\Sigma$. Recall, that in both geometries $N$ is proportional to
$\gamma$ and therefore can be neglected in (\ref{N}):
$$
\dot{\Gamma}(s)=\gamma\wedge\nabla_v\dot{\gamma}=\gamma\wedge
(k\cdot n),
$$
where $k$ is the geodesic curvature of $\gamma$ at the point $r$.

In order to prove Theorem \ref{pm}, we need to prove the equality of
the values of the function $\Psi$ in two points:
$$
P_{-}=\Gamma(s)-\varepsilon \dot{\Gamma}(s),\
P_+=\Gamma(s)+\varepsilon \dot{\Gamma}(s).
$$
We have:
$$
P_-=M-\gamma\wedge(\varepsilon kn)=r\wedge v-r\wedge (\varepsilon
kn)=r\wedge(v-\varepsilon kn),
$$
$$
P_+=M+\gamma\wedge(\varepsilon kn)=r\wedge v+r\wedge (\varepsilon
kn)=r\wedge(v+\varepsilon kn).
$$
The points $P_-,P_+$ do not lie on $\Sigma$ but we can present them
as follows:
$$
P_-=r\wedge(v-\varepsilon kn)=(1+\varepsilon^2k^2)^{1/2}(r\wedge
v_-),
$$
$$
P_+=r\wedge(v+\varepsilon kn)=(1+\varepsilon^2k^2)^{1/2}(r\wedge
v_+).
$$
Here we introduced
$$v_{\pm}:=(1+\varepsilon^2k^2)^{-1/2}(v\pm\varepsilon kn),$$
the two unit vectors in $T_r\Sigma.$ Notice, that the vectors
$v_{\pm}$ obviously obey Birkhoff billiard reflection law. Therefore
we have:
$$\Psi(P_-)=\Psi((1+\varepsilon^2k^2)^{1/2}(r\wedge
v_-))=(1+\varepsilon^2k^2)^{n/2}\Psi(r\wedge
v_-)=(1+\varepsilon^2k^2)^{n/2}\Phi(v_-),$$
$$\Psi(P_+)=\Psi((1+\varepsilon^2k^2)^{1/2}(r\wedge
v_+))=(1+\varepsilon^2k^2)^{n/2}\Psi(r\wedge
v_+)=(1+\varepsilon^2k^2)^{n/2}\Phi(v_+),$$ where we used
homogeneity of $\Psi.$
 Therefore, the claim
follows because by the assumptions $\Phi$ is an integral of Birkhoff
billiard, so that $\Phi(v_-)=\Phi(v_+)$. This proves Theorem
\ref{pm}.
\begin{remark}The following fact is the corollary of the proof. Let $p$ denotes the radial
projection onto the unite Sphere for $K=1$ and onto the de Sitter
surface for $K=-1$. Then we have $$p(P_-):=M_-=r\wedge v_-;\quad
p(P_+):=M_+=r\wedge v_+.$$ It then follows that the geodesic segment
$[M_-;M_+]$ contains $M$ as the middle point (in the sense of
induced metric on the de SItter surface), which is equivalent to the
fact that $v_-$ and $v_+$ respect the Birkhoff billiard reflection
law.
\end{remark}
\section{Proof of Theorem \ref{algebra}.}
\label{sphere} In this section we shall treat the equation
(\ref{pme}) in the both cases of Sphere and Hyperboloid and derive
the equation (\ref{Hess}).

Let us denote by $F\in \mathbb{C}[x_1,x_2,x_3]$ the irreducible
homogenous
 polynomial defining the component
$\tilde{\Gamma}\subset\mathbb{C}P^2$, and let $d=\textrm{deg} F,
d\geq 2$. Since $\Psi$ vanishes on $\Gamma$ we can write
$$\Psi=F^kQ,$$for some integer $k\geq 1$,
so that polynomial $Q$ does not vanish on $\Gamma$ identically.
Notice that $F,Q$ can be assumed to be real since $\Gamma$ is real
curve. If needed we consider a smaller arc where $DF$ does not
vanish and $Q>0.$ We set
$$G=F Q^\frac{1}{k}.$$
Let us remark that $G$ is not a polynomial anymore, but homogeneous
function of degree
\begin{equation}\label{p}p=\frac{n}{k}=d+\frac{1}{k}\textrm{deg} Q\geq 2.
\end{equation} Since $\Psi$ satisfies (\ref{pme}) then also $G$ does.
\subsection{The case of the Sphere.}In this case the tangent vector
$w\in T_r\Gamma, r=(x_1,x_2,x_3)$ can be taken
with the components:
\begin{equation}\label{components}
w_1=x_2 G_{x_3}-x_3G_{x_2};\ w_2=x_3G_{x_1}-x_1G_{x_3};\
w_3=x_1G_{x_2}-x_2G_{x_1}.
\end{equation} In what follows, we pass
from homogeneous functions $G,F,Q,...$ to $g,f,q,...$ defined by
$$f(x,y)=F(x,y,1),\ g(x,y)=G(x,y,1),\ q(x,y)=Q(x,y,1),...,
$$ so that
$$g=fq^{\frac{1}{k}}.$$

Here the mapping $$ (x_1,x_2,x_3)\mapsto (x,y),\ x=\frac{x_1}{x_3},\
y=\frac{x_2}{x_3},\ x_3=(x^2+y^2+1)^{-\frac{1}{2}}
$$
is the central projection of the sphere to the plane $\{z=1\}.$

The derivatives of the functions $G$ and $g$ are related in a usual
way:
\begin{equation}\label{derivative1}
G_{x_1}=x_3^{p-1}g_x=(1+x^2+y^2)^{\frac{1-p}{2}}g_x;
G_{x_2}=x_3^{p-1}g_y=(1+x^2+y^2)^{\frac{1-p}{2}}g_y,
\end{equation}
Using Euler formula for $G$, for all $(x_1,x_2,x_3)\in\Gamma$ so
that $G(x_1,x_2,x_3)=0$, we get also:
\begin{equation}\label{derivative2}
G_{x_3}=-x_3^{p-1}(xg_x+yg_y)=-(1+x^2+y^2)^{\frac{1-p}{2}}(xg_x+yg_y).
\end{equation}
Using (\ref{derivative1}) and (\ref{derivative2}) the components
(\ref{components}) of the vector $w$ take the form:
\begin{equation}\label{w1}
\frac{w_1}{x_3}=-(1+x^2+y^2)^{\frac{1-p}{2}}((xy)g_x+(1+y^2)g_y),
\end{equation}
\begin{equation}\label{w2}
\frac{w_2}{x_3}=(1+x^2+y^2)^{\frac{1-p}{2}}((1+x^2)g_x+(xy)g_y),
\end{equation}
\begin{equation}\label{w3}
\frac{w_3}{x_3}=(1+x^2+y^2)^{\frac{1-p}{2}}(xg_y-yg_x).
\end{equation}

Let us write equation (\ref{pme}) in terms of the function $g$:
\begin{equation}\label{pme1}
(x_3-\varepsilon w_3)^p g\left( \frac{x_1-\varepsilon
w_1}{x_3-\varepsilon w_3},\frac{x_2-\varepsilon w_2}{x_3-\varepsilon
w_3}\right) =
\end{equation}
\begin{equation*}= (x_3+\varepsilon w_3)^p g\left( \frac{x_1+\varepsilon
w_1}{x_3+\varepsilon w_3},\frac{x_2+\varepsilon w_2}{x_3+\varepsilon
w_3}\right).
\end{equation*}
Using formulas (\ref{w1}),(\ref{w2}),(\ref{w3}) equation
(\ref{pme1}) yields:
\begin{equation}\label{pme2}
(1-\mu(xg_y-yg_x))^p g\left(
\frac{x+\mu((xy)g_x+(1+y^2)g_y)}{1-\mu(xg_y-yg_x)},
\frac{y-\mu((1+x^2)g_x+(xy)g_y)}{1-\mu(xg_y-yg_x)}\right)=
\end{equation}
\begin{equation*}=
(1+\mu(xg_y-yg_x))^p g\left(
\frac{x-\mu((xy)g_x+(1+y^2)g_y)}{1+\mu(xg_y-yg_x)},
\frac{y+\mu((1+x^2)g_x+(xy)g_y)}{1+\mu(xg_y-yg_x)}\right),
\end{equation*}
where we denoted by $\mu=\varepsilon(1+x^2+y^2)^{\frac{1-p}{2}}$.
Our next step is to collect and to equate to zero the terms of order
$\mu^3$, neglecting the terms containing $g$, since $g=0$ on the
curve. We have:
\begin{equation}\label{terms}
(1+x^2+y^2)(g_{xxx}g_y^3-3g_{xxy}g_y^2g_x+3g_{xyy}g_yg_x^2-g_{yyy}g_x^3)+
\end{equation}
$$
+3(2-p)(g_{xx}g_y^2-2g_{xy}g_xg_y+g_{yy}g_x^2)(xg_y-yg_x)=0
$$
Denote by $u=g_y\partial_x-g_x\partial_y$ the vector field tangent
to $\{ f=0\}.$ Introduce the quantity
$$
H(g)=g_{xx}g_y^2-2g_{xy}g_xg_y+g_{yy}g_x^2,
$$
and notice that the following two identities hold
$$L_u(x^2+y^2+1)=2(xg_y-yg_x),$$
$$
L_uH(g)=g_{xxx}g_y^3-3g_{xxy}g_y^2g_x+3g_{xyy}g_yg_x^2-g_{yyy}g_x^3.
$$
Using these identities equation (\ref{terms}) can be rewritten as
follows:
\begin{equation}\label{e9}
(x^2+y^2+1)L_u H(g)+\frac{3}{2}(2-p)H(g)L_u(x^2+y^2+1)=0.
\end{equation}
Multiplying (\ref{e9}) by $(x^2+y^2+1)^{\frac{4-3p}{2}}$ we get
\begin{equation}
\label{e10} (x^2+y^2+1)^{\frac{6-3p}{2}}L_uH(g)+
\end{equation}
$$
+\frac{3}{2}(2-p)H(g)(x^2+y^2+1)^{\frac{4-3p}{2}}L_u(x^2+y^2+1)=0.
$$
But the left hand side of (\ref{e10}) is the complete derivative and
thus
$$L_u\left(H(g)(x^2+y^2+1)^{\frac{6-3p}{2}}\right)=0.$$

Since $u$ is a tangent vector field to $\{f=0\}$ then the function
$H(g)(x^2+y^2+1)^{\frac{6-3p}{2}}$ must be a constant on $\{f=0\}$
i.e.
\begin{equation}\label{eq2}
 H(g(x,y))=c_1(x^2+y^2+1)^{\frac{3p-6}{2}},\qquad \forall(x,y)\in\{f=0\},
\end{equation} for some constant $c_1.$
A direct calculation gives that for any function $r$
\begin{equation}\label{cube}
 H(f(x,y)r(x,y))=r^3(x,y)H(f(x,y)), \qquad \forall(x,y)\in\{f=0\},
\end{equation}
and therefore we conclude with the formula
$$
q^\frac{3}{k}(x,y)H(f(x,y))=H(g)=c_1(x^2+y^2+1)^{\frac{3p-6}{2}},\qquad
\forall(x,y)\in\{f=0\}.
$$
Raising to the power $k$ we get
$$
q^3(x,y)H(f(x,y))^k=c_1^k(x^2+y^2+1)^{k\frac{3p-6}{2}},\qquad
\forall(x,y)\in\{f=0\}.
$$
Therefore the difference is a polynomial divisible by $f$.
\begin{equation}\label{d}
q^3(x,y)H(f(x,y))^k-c_1^k(x^2+y^2+1)^{k\frac{3p-6}{2}}=f\cdot
r_1,\qquad \forall(x,y)\in\mathbb{C}^2, \end{equation} for some
polynomial $r_1.$ Moreover, using (\ref{d})one can compute the
degrees of the terms at the left hand side of (\ref{d}) and to get
homogeneous version of (\ref{d}):\begin{equation}\label{e6}
Q^3(x,y,z)H(F(x,y,z))^k-c_1^k(x^2+y^2+z^2)^{k\frac{3p-6}{2}}z^{2k}=F\cdot
R_1,
\end{equation}
which is valid for all $(x,y,z)\in\mathbb{C}^3.$ Now we use the
identities (see \cite{W}):
$$
 {\rm Hess}(F)=\frac{(d-1)^2}{z^2}
  \det\left(
  \begin{array}{ccc}
    F_{xx} & F_{xy} & F_{x} \\
    F_{xy} & F_{yy} & F_{y} \\
    F_{x} & F_{y} & \frac{d}{d-1}F\\
  \end{array}
\right)=$$
\begin{equation}\label{Hf}
=\frac{(d-1)^2}{z^2}\left(\frac{d}{d-1}F
 (F_{xx}F_{yy}-F_{xy}^2)-H(F)\right).
\end{equation}
Notice that in (\ref{Hf}) and in (\ref{e6}) $H(F(x,y,z))$  is
computed with $z$ being a parameter. The last step is just to
substitute the expression for $H(F)$ via ${\rm Hess}(F)$ from
(\ref{Hf}) into (\ref{e6}) in order to get:
\begin{equation}\label{Hess1}
Q^3{\rm (Hess (F))}^k-c(x^2+y^2+z^2)^{k\frac{3p-6}{2}}=F\cdot
R,\quad \forall(x,y,z)\in\mathbb{C}^3.
\end{equation}
But this is exactly what it is claimed in
 (\ref{Hess}).
\subsection{The case of Hyperboloid}
\label{hyp} We proceed exactly as in previous Subsection with
obvious modifications. We keep the same notation as in the beginning
of the Section.

Recall that in the Hyperbolic case $\Sigma$ is the upper sheet of
Hyperboloid, the curve $\gamma\subset\Sigma$ and the curve $\Gamma$
lies on the one sheeted Hyperboloid
$$x_1^2+x_2^2-x_3^2=1.$$
Therefore, the tangent vector $w\in T_r\Gamma, r=(x_1,x_2,x_3)$ can
be taken with the components:
\begin{equation}\label{hypcomponents}
w_1=x_2 G_{x_3}+x_3G_{x_2};\ w_2=-x_3G_{x_1}-x_1G_{x_3};\
w_3=x_1G_{x_2}-x_2G_{x_1}.\end{equation}

As before we pass from homogeneous functions $G,F,Q$ to $g,f,q$
defined above:
$$f(x,y)=F(x,y,1),\ g(x,y)=G(x,y,1),\ q(x,y)=Q(x,y,1),..,
$$
where the mapping $$ (x_1,x_2,x_3)\mapsto (x,y),\
x=\frac{x_1}{x_3},\ y=\frac{x_2}{x_3}
$$ and $x_3=(x^2+y^2-1)^{-\frac{1}{2}}$
is the central projection of the one sheeted hyperboloid to the
plane $\{z=1\}.$ The derivatives of the functions $G$ and $g$ are
related in this case as follows:
\begin{equation}\label{hypderivative1}
G_{x_1}=x_3^{p-1}g_x=(x^2+y^2-1)^{\frac{1-p}{2}}g_x,
G_{x_2}=x_3^{p-1}g_y=(x^2+y^2-1)^{\frac{1-p}{2}}g_y.
\end{equation}
Using Euler formula for $G$, for all $(x_1,x_2,x_3)\in\Gamma$ so
that $G(x_1,x_2,x_3)=0$, we get also:
\begin{equation}\label{hypderivative2}
G_{x_3}=-x_3^{p-1}(xg_x+yg_y)=-(x^2+y^2-1)^{\frac{1-p}{2}}(xg_x+yg_y).
\end{equation}
Using (\ref{hypderivative1}) and (\ref{hypderivative2}) the
components (\ref{hypcomponents}) of the vector $w$ in this case take
the form:
\begin{equation}\label{hypw1}
\frac{w_1}{x_3}=-(x^2+y^2-1)^{\frac{1-p}{2}}((xy)g_x+(y^2-1)g_y),
\end{equation}
\begin{equation}\label{hypw2}
\frac{w_2}{x_3}=(x^2+y^2-1)^{\frac{1-p}{2}}((x^2-1)g_x+(xy)g_y),
\end{equation}
\begin{equation}\label{hypw3}
\frac{w_3}{x_3}=(x^2+y^2-1)^{\frac{1-p}{2}}(xg_y-yg_x).
\end{equation}

Therefore equation (\ref{pme1}) in the case of Hyperboloid can be
presented using formulas (\ref{hypw1}), (\ref{hypw2}),
(\ref{hypw3}):
\begin{equation}\label{hyppme2}
(1-\mu(xg_y-yg_x))^p g\left(
\frac{x+\mu((xy)g_x+(y^2-1)g_y)}{1-\mu(xg_y-yg_x)},
\frac{y-\mu((x^2-1)g_x+(xy)g_y)}{1-\mu(xg_y-yg_x)}\right)=
\end{equation}
\begin{equation*}=
(1+\mu(xg_y-yg_x))^p g\left(
\frac{x-\mu((xy)g_x+(y^2-1)g_y)}{1+\mu(xg_y-yg_x)},
\frac{y+\mu((x^2-1)g_x+(xy)g_y)}{1+\mu(xg_y-yg_x)}\right),
\end{equation*}
where we denoted $\mu=\varepsilon(x^2+y^2-1)^{\frac{1-p}{2}}$. Our
next step is to collect and equate to zero the terms of order
$\mu^3$, neglecting the terms containing $g$:
\begin{equation}\label{hypterms}
(x^2+y^2-1)(g_{xxx}g_y^3-3g_{xxy}g_y^2g_x+3g_{xyy}g_y
g_x^2-g_{yyy}g_x^3)+
\end{equation}
$$
+3(2-p)(g_{xx}g_y^2-2g_{xy}g_xg_y+g_{yy}g_x^2)(xg_y-yg_x)=0.
$$

As before we denote by $u=g_y\partial_x-g_x\partial_y$ be the vector
field tangent to $\{ f=0\}$. Recall the notation,
$$
H(g)=g_{xx}g_y^2-2g_{xy}g_xg_y+g_{yy}g_x^2,
$$
and use the two identities:
$$L_u(x^2+y^2-1)=2(xg_y-yg_x),$$
$$
L_uH(g)=g_{xxx}g_y^3-3g_{xxy}g_y^2g_x+3g_{xyy}g_y
g_x^2-g_{yyy}g_x^3.
$$
Using these identities equation (\ref{hypterms}) can be rewritten as
follows:
\begin{equation}\label{hype9}
(x^2+y^2-1)L_u H(g)+\frac{3}{2}(2-p)H(g)L_u(x^2+y^2-1)=0.
\end{equation}
Multiplying again (\ref{hype9}) by $(x^2+y^2-1)^{\frac{4-3p}{2}}$ we
get
\begin{equation}
\label{hype10} (x^2+y^2-1)^{\frac{6-3p}{2}}L_uH(g)+
\end{equation}
$$
+\frac{3}{2}(2-p)H(g)(x^2+y^2-1)^{\frac{4-3p}{2}}L_u(x^2+y^2-1)=0.
$$
But the left hand side of (\ref{hype10}) is the complete derivative
and thus
$$L_u\left(H(g)(x^2+y^2-1)^{\frac{6-3p}{2}}\right)=0.$$
Since $u$ is a tangent vector field to $\{f=0\}$, then the function
\newline $H(g)(x^2+y^2-1)^{\frac{6-3p}{2}}$ must be a constant on
$\{f=0\}$, i.e.
\begin{equation}\label{hypeq2}
 H(g(x,y))=c_1(x^2+y^2-1)^{\frac{3p-6}{2}},\qquad \forall(x,y)\in\{f=0\},
\end{equation} for some constant $c_1.$ Using again the identity
(\ref{cube}) we come to:
$$
q^\frac{3}{k}(x,y)H(f(x,y))=c_1(x^2+y^2-1)^{\frac{3p-6}{2}},\qquad
\forall(x,y)\in\{f=0\}.
$$
Raising to the power $k$ we get
$$
q^3(x,y)H(f(x,y))^k=c_1^k(x^2+y^2-1)^{k\frac{3p-6}{2}},\qquad
\forall(x,y)\in\{f=0\}.
$$
Therefore the difference is a polynomial divisible by $f$.
$$
q^3(x,y)H(f(x,y))^k-c_1^k(x^2+y^2-1)^{k\frac{3p-6}{2}}=f\cdot
r_1,\qquad \forall(x,y)\in\mathbb{C}^2,
$$ for some polynomial $r_1.$
One can compute the degrees using (\ref{d}) and to get  homogeneous
version of the last identity:
\begin{equation}\label{hype6}
Q^3(x,y,z)H(F(x,y,z))^k-c_1^k(x^2+y^2-z^2)^{k\frac{3p-6}{2}}z^{2k}=F\cdot
R_1,
\end{equation}
which is valid for all $(x,y,z)\in\mathbb{C}^3.$ Using again the
identities (\ref{Hf}) we substitute the expression for
$H(F)=\frac{z^2}{(d-1)^2} \rm Hess(F)$ which is valid on $\{F=0\}$
into (\ref{hype6}). We get:
\begin{equation}\label{hypHess1}
Q^3{\rm (Hess (F))}^k-c(x^2+y^2-z^2)^{k\frac{3p-6}{2}}=F\cdot
R,\quad \forall(x,y,z)\in\mathbb{C}^3.
\end{equation}
But this is exactly the equation of
 (\ref{Hess}) for the Hyperbolic case. This completes the proof of Theorem \ref{algebra} in the Hyperbolic case.
\begin{remark}It would be very natural to collect terms of
$\varepsilon^3$ of the equation (\ref{pme}) written in homogeneous
coordinates without going to affine chart. One would expect that
using Euler formula for $F$ and maybe for its derivatives one would
be able to get directly to the equation (\ref{Hess}). Surprisingly,
this attempt leads to very heavy calculations, which we were not
able to perform.
\end{remark}

\end{document}